\documentclass[12pt]{article}
\usepackage{latexsym}
\usepackage{amsmath}
\usepackage{amssymb}
\usepackage{amsfonts}

\numberwithin{equation}{section}

\newtheorem{theorem}{Theorem}[section]
\newtheorem{lemma}[theorem]{Lemma}

\newtheorem{cor}[theorem]{Corollary}
\newtheorem{remark}[theorem]{Remark}
\newcommand{\eproof}{{\mbox{\ }~\hfill
\mbox{\large $\Box$} \par \vskip 10pt}}

\newcommand{\R}{{\mathbb R}}

\newcommand{\pf}{\noindent{\bf Proof}}
\renewcommand{\div}{{\rm div}}

\title{Optimal three-ball inequalities and quantitative uniqueness for the Stokes system}
\author{Ching-Lung Lin\thanks{Department of Mathematics, National Cheng Kung University,
Tainan 701, Taiwan. Email:cllin2@mail.ncku.edu.tw}\qquad Jenn-Nan
Wang\thanks{Department of Mathematics, Taida Institute of
Mathematical Sciences, NCTS (Taipei), National Taiwan University,
Taipei 106, Taiwan. Email: jnwang@math.ntu.edu.tw}}

\date{}

\begin{document}
\maketitle

\begin{abstract}
In this paper we study the local behavior of a solution to the
Stokes system with singular coefficients. One of the main results is
the bound on the vanishing order of a nontrivial solution to the
Stokes system, which is a quantitative version of the strong unique
continuation property. Our proof relies on some delicate
Carleman-type estimates. We first use these estimates to derive
crucial \emph{optimal} three-ball inequalities. Taking advantage of
the optimality, we then derive an upper bound on the vanishing order
of any nontrivial solution to the Stokes system from those
three-ball inequalities.
\end{abstract}

\section{Introduction}\label{sec1}
\setcounter{equation}{0}

Assume that $\Omega$ is a connected open set containing $0$ in
$\R^n$ with $n\geq 2$. In this paper we are interested in the local
behavior of $(u,p)$ satisfying the following Stokes system:
\begin{equation}\label{1.1}
\begin{cases}
\begin{array}{l}
\Delta u+A(x)\cdot \nabla u+\nabla p=0\\
\div u=0,
\end{array}
\end{cases}
\end{equation}
where $A$ is measurable satisfying
\begin{equation}\label{asing}
|A(x)|\leq C_0|\log |x||^{-3}|x|^{-1}\quad\forall\ x\in\Omega
\end{equation}
and $A\cdot \nabla u=(A\cdot \nabla u_1,\cdots,A\cdot \nabla u_n)$.

For the Stokes system \eqref{1.1} with essentially bounded
coefficients $A(x)$, the weak unique continuation property has been
shown by Fabre and Lebeau \cite{fale}. On the other hand, when
$A(x)$ satisfies $|A(x)|=O(|x|^{-1+\epsilon})$ with $\epsilon>0$,
the strong unique continuation property was proved by Regbaoui
\cite{Reg2}. The results in \cite{fale} and \cite{Reg2} concern only
the qualitative unique continuation theorem. In this work we aim to
derive a quantitative estimate of the strong unique continuation for
\eqref{1.1}.

For the second order elliptic operator, using Carleman or frequency
functions methods, quantitative estimates of the strong unique
continuation (in the form of doubling inequality) under different
assumptions on coefficients were derived in \cite{dofe}, \cite{gl1},
\cite{gl2}, \cite{lin1}, \cite{lin2}. For the power of Laplacian, a
quantitative estimate was obtained in \cite{lin3}. We refer to
\cite{lin2} and references therein for the development of this
investigation.

Since there is no equation for $p$ in the Stokes system \eqref{1.1},
to prove the unique continuation theorem for \eqref{1.1}, one
usually apply the divergence on the first equation and obtain
\begin{equation}\label{eq0}
\Delta p+\div (A(x)\cdot\nabla u)=0.
\end{equation}
However, the first equation of \eqref{1.1} and \eqref{eq0} do not
give us a decoupled system. The frequency functions method does not
seem to work in this case. So we prove our results along the line of
Carleman's method. On the other hand, since the coefficient $A(x)$
is more singular than the one considered in \cite{Reg2}.
Carleman-type estimates derived in \cite{Reg2} can not be applied to
the case here. Hence we need to derive new Carleman-type estimates
for our purpose. The key is to use weights which are slightly less
singular than the negative powers of $|x|$ (see estimates
\eqref{2.4} and \eqref{2.14}). The estimate \eqref{2.14} is to
handle \eqref{eq0} and the idea is due to Fabre and Lebeau
\cite{fale}. It is tempting to derive doubling inequalities for
\eqref{1.1} by \eqref{2.4} and \eqref{2.14} using the ideas in
\cite{lin2} or \cite{lin3}. But this seems hard to reach with
estimates \eqref{2.4}, \eqref{2.14}. One of the difficulties is the
appearance of the parameter $\beta$ on the right hand side of
\eqref{2.14}.

Even though we are not able to prove doubling inequalities for
\eqref{1.1}, we can derive certain three-ball inequalities which are
\emph{optimal} in the sense explained in \cite{efv} using
\eqref{2.4} and \eqref{2.14}. We would like to remark that usually
the three-ball inequality can be regarded as the quantitative
estimate of the weak unique continuation property. However, when the
three-ball inequality is optimal, one is able to deduce the strong
unique continuation from it. It seems reasonable to expect that one
could derive a bound on the vanishing order of a nontrivial solution
from the optimal three-ball inequality. A recent result by Bourgain
and Kenig \cite{bou} (more precisely, Kenig's lecture notes for 2006
CNA Summer School \cite{kenigcna}) indicates that this is indeed
possible, at least for the Schr\"odinger operator. In this paper, we
show that by the optimal three-ball inequality, we can obtain a
bound on the vanishing order of a nontrivial solution to \eqref{1.1}
containing "nearly" optimal singular coefficients. Finally, we would
like to mention that quantitative estimates of the strong unique
continuation are useful in studying the nodal sets of solutions for
elliptic or parabolic equations \cite{hasi}, \cite{linf2}, or the
inverse problem \cite{abrv}.

We now state main results of this paper. Their proofs will be given
in the subsequent sections. Assume that there exists $0<R_0\le 1$
such that $B_{R_0}\subset\Omega$. Hereafter $B_r$ denotes an open
ball of radius $r>0$ centered at the origin. Also, we let
$U(x)=[|x|^{4}|\nabla p|^2+|x|^{2}|p|^2+|u|^2]^{1/2}$.
\begin{theorem}\label{thm1.1}
There exists a positive number $\tilde{R}<1$, depending only on $n$,
such that if $\ 0<R_1<R_2<R_3\leq R_0$ and
$R_1/R_3<R_2/R_3<\tilde{R}$, then
\begin{equation}\label{1.2}
\int_{|x|<R_2}|U|^2dx\leq
{C}\left(\int_{|x|<R_1}|U|^2dx\right)^{\tau}\left(\int_{|x|<{R_3}}|U|^2dx\right)^{1-\tau}
\end{equation}
for $(u,p)\in (H^1({B}_{R_0}))^{n+1}$ satisfying \eqref{1.1} in
${B}_{R_0}$, where the constant ${C}$ depends on $R_2/R_3$, $n$, and
$0<\tau<1$ depends on $R_1/R_3$, $R_2/R_3$, $n$. Moreover, for fixed
$R_2$ and $R_3$, the exponent $\tau$ behaves like $1/(-\log R_1)$
when $R_1$ is sufficiently small.
\end{theorem}
\begin{remark}\label{rem1.1}
It is important to emphasize that $C$ is independent of $R_1$ and
$\tau$ has the asymptotic $(-\log R_1)^{-1}$. These facts are
crucial in deriving an vanishing order of a nontrivial $(u,p)$ to
\eqref{1.1}. Due to the behavior of $\tau$, the three-ball
inequality is called optimal {\rm\cite{efv}}.
\end{remark}
\begin{remark}
We want to say a few words about the appearance of $\nabla p$ term
in $U$. In the derivation of the three-ball inequality \eqref{1.2},
it is crucial to control $\nabla u$ in a smaller region by
quantities of $u$ and $p$ in a bigger region {\rm(}see
\eqref{3.1}{\rm )}. Roughly speaking, this is an interior estimate
for $\nabla u$. In view of the first equation of \eqref{1.1},
$\nabla p$ needs be included in this estimate.
\end{remark}

\begin{theorem}\label{thm1.2}
Let $(u,p)\in (H^1({B}_{R_0}))^{n+1}$ be a nontrivial solution to
\eqref{1.1}, i.e, $(u,p)\ne(0,0)$, then there exist positive
constants $K$ and $m$, depending on $n$ and $(u,p)$, such that
\begin{equation}\label{1.3}
\int_{|x|<R}|U|^2 dx\ge KR^m
\end{equation}
for all $R$ sufficiently small.
\end{theorem}
\begin{remark}\label{rem1.2}
Based on Theorem~\ref{thm1.1}, the constants $K$ and $m$ in
\eqref{1.3} are given by
$$
K=\int_{|x|<R_3}|U|^2dx
$$
and
$$
m=\tilde
C\log\Big{(}\frac{\int_{|x|<R_3}|U|^2dx}{\int_{|x|<R_2}|U|^2dx}\Big{)},
$$
where $\tilde C$ is a positive constant depending on $n$ and
$R_2/R_3$.
\end{remark}
\begin{cor}\label{cor1}
Let $(u,p)\in (H^1_{loc}(\Omega))^{n}\times L_{loc}^2(\Omega)$ be a
solution of \eqref{1.1} with $A$ satisfying \eqref{asing}. Assume
that $(u,p)$ vanishes of infinite order at the origin, i.e., for all
$N>0$,
\begin{equation}\label{info}
\int_{|x|<R}(|u|^2+|p|^2)dx=O(R^N)\quad\text{as}\quad R\to 0.
\end{equation}
Then $(u,p)\equiv 0$ in $\Omega$.
\end{cor}
This corollary is a small improvement of the strong unique
continuation property for the Stokes system proved in \cite{Reg2}
where $|A(x)|=O(|x|^{-1+\epsilon})$ with $\epsilon>0$.

This paper is organized as follows. In Section~2, we derive suitable
Carleman-type estimates. A technical interior estimate is proved in
Section~3. Section~4 is devoted to the proofs of
Theorem~\ref{thm1.1}, \ref{thm1.2}, and Corollary~\ref{cor1}.

\section{Carleman estimates}\label{sec2}
\setcounter{equation}{0}

Similar to the arguments used in \cite{Hor1}, we introduce polar
coordinates in ${\mathbb R}^n \backslash {\{0\}}$ by setting $x=r
\omega$, with $r=|x|$, $\omega=(\omega_1,\cdots,\omega_n)\in
S^{n-1}$. Furthermore, using new coordinate $t=\log r$, we can see
that
$$\frac{\partial}{\partial x_j}=e^{-t}(\omega_j \partial_t +\Omega_j),\quad 1\le j\le n,$$
where $\Omega_j$ is a vector field in $S^{n-1}$. We could check that
the vector fields $\Omega_j$ satisfy
$$\sum_{j=1}^n\omega_j\Omega_j=0\quad\text{and}\quad\sum_{j=1}^n\Omega_j\omega_j=n-1.$$
Since $r\rightarrow 0$ iff $t\rightarrow {-\infty}$, we are mainly
interested in values of $t$ near $-\infty$.

It is easy to see that
\begin{equation*}
\frac{\partial ^2}{\partial x_j \partial x_{\ell}}=e^{-2t}(\omega_j
\partial_t -\omega_j +\Omega_j)(\omega_{\ell} \partial_t +\Omega_{\ell}),\quad 1\le j,\ell\le n.
\end{equation*}
and, therefore, the Laplacian becomes
\begin{equation}\label{2.1}
e^{2t}\Delta =\partial^2_t +(n-2)\partial_t +\Delta_\omega,
\end{equation}
where $\Delta_\omega=\Sigma^n_{j=1}\Omega^2_j$ denotes the
Laplace-Beltrami operator on $S^{n-1}$. We recall that the
eigenvalues of $-\Delta_\omega$ are $k(k+n-2), k\in \mathbb{N}$, and
the corresponding eigenspaces are $E_k$, where $E_k$ is the space of
spherical harmonics of degree $k$. It follows that
\begin{equation}\label{2.2}
\iint |\Delta_\omega v|^2 dt d\omega=\sum_{k\geq 0}k^2(k+n-2)^2
\iint | v_k |^2 dt d\omega
\end{equation}
and
\begin{equation}\label{2.3}
\sum_{j}\iint |\Omega_j v|^2 dt d\omega=\sum_{k\geq 0}k(k+n-2)\iint|
v_k |^2 dt d\omega,
\end{equation}
where $v_k$ is the projection of $v$ onto $E_k$. Let
$$
\Lambda=\sqrt{\frac{(n-2)^2}{4}-\Delta_{\omega}},
$$
then $\Lambda$ is an elliptic first-order positive
pseudodifferential operator in $L^2(S^{n-1})$. The eigenvalues of
$\Lambda$ are $k+\frac{n-2}{2}$ and the corresponding eigenspaces
are $E_k$. Denote
$$
L^{\pm}=\partial_t+\frac{n-2}{2}\pm\Lambda.
$$
Then it follows from \eqref{2.1} that
\begin{equation*}
e^{2t}\Delta =L^+L^-=L^-L^+.
\end{equation*}

Motivated by the ideas in \cite{Reg}, we will derive Carleman-type
estimates with weights $\varphi_{\beta}=\varphi_{\beta}(x) =\exp
(-\beta\tilde{\psi}(x))$, where $\beta>0$ and $\tilde{\psi}(x)=\log
|x|+\log((\log |x|)^2)$. Note that $\varphi_{\beta}$ is less
singular than $|x|^{-\beta}$, For simplicity, we denote
$\psi(t)=t+\log t^2$, i.e., $\tilde{\psi}(x)=\psi(\log|x|)$. From
now on, the notation $X\lesssim Y$ or $X\gtrsim Y$ means that $X\le
CY$ or $X\ge CY$ with some constant $C$ depending only on $n$.

\begin{lemma}\label{lem2.1}
There exist a sufficiently small $r_0>0$ depending on $n$ and a
sufficiently large $\beta_0>1$ depending on $n$ such that for all
$u\in U_{r_0}$ and $\beta\geq \beta_0$, we have that
\begin{equation}\label{2.4}
\beta  \int \varphi^2_\beta (\log|x|)^{-2}|x|^{-n}(|x|^{2}|\nabla
u|^2+|u|^2)dx\lesssim\int \varphi^2_\beta|x|^{-n}|x|^{4}|\Delta u|^2
dx,
\end{equation}
where $U_{r_0}=\{u\in C_0^{\infty}(\R^n\setminus\{0\}): \mbox{\rm
supp}(u)\subset B_{r_0}\}$.
\end{lemma}

\pf. By the polar coordinate system described above, we have
\begin{eqnarray}\label{2.5}
&&\int \varphi^2_\beta |x|^{4-n}|\Delta u|^2 dx\notag\\
&=&\iint e^{-2\beta\psi(t)}  e^{4t}|\Delta u|^2 dt d\omega\notag\\
&=&\iint| e^{-\beta\psi(t)} e^{2t} \Delta u|^2 dt d\omega.
\end{eqnarray}
If we set $u=e^{\beta\psi(t)}v$ and use \eqref{2.1}, then
\begin{equation}\label{2.6}
e^{-\beta\psi(t)}e^{2t} \Delta
u=\partial_t^2v+b\partial_tv+av+\Delta_\omega v=:P_\beta v,
\end{equation}
where
$a=(1+2t^{-1})^2\beta^2+(n-2)\beta+2(n-2)t^{-1}\beta-2t^{-2}\beta$
and $b=n-2+2\beta+4t^{-1}\beta$. By \eqref{2.5} and \eqref{2.6},
\eqref{2.4} holds if for $t$ near $-\infty$ we have
\begin{equation}\label{2.7}
\sum_{j+|\alpha|\leq1}\beta^{3-2|\alpha|}\iint
|t|^{-2}|\partial_t^j\Omega^{\alpha}v|^2 dt d\omega\leq \tilde C_1 \iint
|P_\beta v|^2 dt d\omega,
\end{equation}
where $\tilde C_1$ is a positive constant depending on $n$.

From \eqref{2.6}, using the integration by parts, for $t<t_0$ and
$\beta>\beta_0$, where $t_0<-1$ and $\beta_0>0$ depend on $n$, we
have that
\begin{eqnarray}\label{2.8}
&& \iint |P_\beta v|^2 dt d\omega\notag\\
&=&\iint
|\partial_t^2v|^2dtd\omega+\iint|b\partial_tv|^2dtd\omega+\iint|av|^2dtd\omega+\iint |\Delta_{\omega}v|^2dtd\omega\notag\\
&&-\iint\partial_tb|\partial_tv|^2dtd\omega-2\iint
a|\partial_tv|^2dtd\omega+\iint\partial_t^2a|v|^2dtd\omega\notag\\
&&-\iint\partial_t(ab)|v|^2dtd\omega+2\sum_j\iint|\partial_t\Omega_jv|^2dtd\omega\notag\\
&&+\sum_j\iint\partial_tb|\Omega_jv|^2dtd\omega-2\sum_j\iint a|\Omega_jv|^2dtd\omega\notag\\
&\ge&\iint
|\Delta_{\omega}v|^2dtd\omega+\iint\{b^2-\partial_tb-2a\}|\partial_tv|^2dtd\omega\notag\\
&&+\sum_j\iint\{\partial_tb-2a\}|\Omega_jv|^2dtd\omega+\iint\{a^2+\partial_t^2a-\partial_t(ab)\}|v|^2dtd\omega\notag\\
&\ge&\iint
|\Delta_{\omega}v|^2dtd\omega+\sum_j\iint\{-4t^{-2}\beta-2a\}|\Omega_jv|^2dtd\omega\notag\\
&&+\iint\{a^2+11t^{-2}\beta^3\}|v|^2dtd\omega+\iint\beta^{2}|\partial_tv|^2dtd\omega.
\end{eqnarray}
In view of \eqref{2.8}, using \eqref{2.2},\eqref{2.3}, we see that
\begin{eqnarray}\label{2.88}
&&\iint |\Delta_{\omega}v|^2dtd\omega-2\sum_j\iint
a|\Omega_jv|^2dtd\omega+\iint a^2|v|^2dtd\omega\notag\\
&=&\sum_{k\ge 0}\iint [a-k(k+n-2)]^2|v_k|^2dtd\omega.
\end{eqnarray}

Substituting \eqref{2.88} into \eqref{2.8} yields
\begin{eqnarray}\label{2.99}
&& \iint |P_\beta v|^2 dt d\omega\notag\\
&\ge&\sum_{k\geq 0}\iint \{11t^{-2}\beta^3-4t^{-2}\beta k(k+n-2)+[a-k(k+n-2)]^2\}|v_k|^2dt d\omega\notag\\
&&+\iint\beta^{2}|\partial_tv|^2dtd\omega\notag\\
&=&\big{(}\sum_{k, k(k+n-2)\ge 2\beta^2}+\sum_{k, k(k+n-2)< 2\beta^2}\big{)}\iint \{11t^{-2}\beta^3-4t^{-2}\beta k(k+n-2)\notag\\
&&\hspace{10mm}+[a-k(k+n-2)]^2\}|v_k|^2dt d\omega+\iint\beta^{2}|\partial_tv|^2dtd\omega.
\end{eqnarray}
For $k$ such that $k(k+n-2)<2\beta^2$, we have
\begin{equation}\label{2.100}
11t^{-2}\beta^3-4t^{-2}\beta k(k+n-2)\ge t^{-2}\beta^3+t^{-2}\beta k(k+n-2).
\end{equation}
On the other hand, if $2\beta^2<k(k+n-2)$, then, by taking $t$ even smaller, if necessary, we get that
\begin{equation}\label{2.200}
-4t^{-2}\beta k(k+n-2)+[a-k(k+n-2)]^2\gtrsim t^{-2}\beta k(k+n-2).
\end{equation}
Finally, using formula \eqref{2.3} and estimates \eqref{2.100},
\eqref{2.200} in \eqref{2.99}, we immediately obtain \eqref{2.7} and
the proof of the lemma is complete.\eproof

To handle the auxiliary equation corresponding to the pressure $p$, we need another Carleman estimate. The derivation here follows the line in \cite{Reg2}.
\begin{lemma}\label{lem2.2}
There exists a sufficiently small number $t_0<0$ depending on $n$
such that for all $u\in V_{t_0}$, $\beta> 1$, we have that
\begin{equation}\label{2.9}
\sum_{j+|\alpha|\leq1}\beta^{1-2(j+|\alpha|)}\iint
t^{-2}\varphi^2_\beta|\partial_t^j\Omega^\alpha u|^2 dt d\omega
\lesssim \iint\varphi^2_\beta |L^- u|^2 dtd\omega,
\end{equation}
where $V_{t_0}=\{u(t,\omega)\in C_0^{\infty}((-\infty,t_0)\times
S^{n-1})\}$.
\end{lemma}
\pf.  If we set $u=e^{\beta\psi(t)}v$, then simple integration by parts implies
\begin{eqnarray*}
&&\iint\varphi^2_\beta |L^- u|^2 dtd\omega\notag\\
&=&\iint  |\partial_t v-\Lambda v+\beta v+2\beta t^{-1}v+(n-2)v/2|^2dt d\omega\notag\\
&=&\iint |\partial_t v|^2 dt d\omega+\iint |-\Lambda v+\beta
v+2\beta t^{-1}v+(n-2)v/2|^2 dt d\omega\notag\\
&&+\beta\iint  t^{-2}|v|^2 dtd\omega.
\end{eqnarray*}
By the definition of $\Lambda$, we have
\begin{eqnarray*}
&&\iint |-\Lambda v+\beta
v+2\beta t^{-1}v+(n-2)v/2|^2 dt d\omega\notag\\
&=&\sum_{k\geq 0}\iint  |-k v_k+\beta v_k+2\beta t^{-1}v_k|^2
dtd\omega\notag\\
&=&\sum_{k\geq 0}\iint  (-k +\beta+2\beta t^{-1})^2|v_k|^2
dtd\omega,
\end{eqnarray*}
where, as before, $v_k$ is the projection of $v$ on $E_k$. Note
that
$$
(-k +\beta+2\beta t^{-1})^2+\beta t^{-2}\ge\frac{1}{8\beta}(2\beta
t^{-1})^2+\frac{1}{16\beta}(\beta-k)^2.
$$
Considering $\beta>(1/2)k$ and $\beta\le (1/2)k$, we can get that
\begin{eqnarray}\label{2.12}
&&\iint \varphi^2_\beta |L^- u|^2 dtd\omega\notag\\
&=&\iint|\partial_t v|^2 dt d\omega+\Sigma_{k\geq 0}\iint[(-k
+\beta+2\beta t^{-1})^2+\beta
t^{-2}]|v_k|^2 dtd\omega\notag\\
&\gtrsim& \iint|\partial_t v|^2 dt d\omega+\Sigma_{k\geq 0}\iint
(\beta^{-1}t^{-2}k(k+n-2) +\beta t^{-2})|v_k|^2 dtd\omega.\notag\\
\end{eqnarray}
The estimate \eqref{2.9} then follows from \eqref{2.3}.\eproof

Next we need a technical lemma. We then use this lemma to derive
another Carleman estimate.
\begin{lemma}\label{lem2.3}
There exists a sufficiently small number $t_1<-2$ depending on $n$
such that for all $u\in V_{t_1}$, $g=(g_0,g_1,\cdots,g_n)\in
(V_{t_1})^{n+1}$ and $\beta> 0$, we have that
\begin{equation*}
\iint \varphi^2_\beta |u|^2 dt d\omega \lesssim \iint
\varphi^2_\beta (|L^+
u+\partial_tg_0+\sum_{j=1}^n\Omega_jg_j|^2+\|g\|^2) dtd\omega.
\end{equation*}
\end{lemma}
\pf. This lemma can be proved by exactly the same arguments used in
Lemma 2.2 of \cite{Reg2}. So we omit the proof here.\eproof

\begin{lemma}\label{lem2.4}
There exist a sufficiently small number $r_1>0$ depending on $n$ and
a sufficiently large number $\beta_1>2$ depending on $n$ such that
for all $w\in U_{r_1}$ and $f=(f_1,\cdots,f_n)\in (U_{r_1})^{n}$,
$\beta\geq \beta_1$, we have that
\begin{eqnarray}\label{2.14}
&&\int\varphi^2_\beta (\log|x|)^2(|x|^{4-n}|\nabla w|^2+|x|^{2-n}|w|^2)dx\notag\\
&\lesssim& \beta\int \varphi^2_\beta
(\log|x|)^{4}|x|^{2-n}[(|x|^{2}\Delta w+|x|{\rm div} f)^2+\|f\|^2
]dx,
\end{eqnarray}
where $U_{r_1}$ is defined as in Lemma~\ref{lem2.1}.
\end{lemma}
\pf. Replacing $\beta$ by $\beta+2$ in \eqref{2.14}, we see that it
suffices to prove
\begin{eqnarray}\label{2.15}
&&\int\varphi^2_\beta (\log|x|)^{-2}(|x|^{2}|\nabla w|^2+|w|^2)|x|^{-n}dx\notag\\
&\lesssim& \beta\int \varphi^2_\beta [(|x|^{2}\Delta w+|x|{\rm div}
f)^2+\| f\|^2 ] |x|^{-n}dx.
\end{eqnarray}
Working in polar coordinates and using the relation $e^{2t}\Delta
=L^+L^-$, \eqref{2.15} is equivalent to
\begin{eqnarray}\label{2.16}
&&\sum_{j+|\alpha|\leq1}\iint
t^{-2}\varphi^2_\beta|\partial_t^j\Omega^\alpha u|^2dtd\omega\notag\\
&\lesssim& \beta\iint \varphi^2_\beta (|L^+L^-
w+\partial_t(\sum_{j=1}^n \omega_jf_j)+\sum_{j=1}^n\Omega_jf_j|^2+\|
f\|^2)dtd\omega.
\end{eqnarray}
Applying Lemma \ref{lem2.3} to $u=L^-w$ and $g=(\sum_{j=1}^n
\omega_jf_j,f_1,\cdots,f_n)$ yields
\begin{eqnarray}\label{2.17}
&&\beta\iint
\varphi^2_\beta |L^-w|^2 dt d\omega\notag\\
&\lesssim& \beta\iint\varphi^2_\beta (|L^+
L^-w+\partial_t(\sum_{j=1}^n
\omega_jf_j)+\sum_{j=1}^n\Omega_jf_j|^2+\|f\|^2) dtd\omega.
\end{eqnarray}
Now \eqref{2.16} is an easy consequence of \eqref{2.9} and
\eqref{2.17}.\eproof

\section{Interior estimates}\label{sec3}
\setcounter{equation}{0}

To establish the three-ball inequality for \eqref{1.1}, the
following interior estimate is useful.
\begin{lemma}\label{lem3.1}
Let $(u,p)\in (H^1_{loc}(\Omega))^{n+1}$ be a solution of
\eqref{1.1}. Then for any $0<a_3<a_1<a_2<a_4$ such that
$B_{a_4r}\subset\Omega$ and $|a_4r|<1$, we have
\begin{equation}\label{3.1}
\int_{a_1r<|x|<a_2r}|x|^{2}|\nabla u|^2dx\le
C\int_{a_3r<|x|<a_4r}(|x|^{4}|\nabla p|^2+|u|^2)dx,
\end{equation}
where the constant $C$ is independent of $r$ and $(u,p)$.
\end{lemma}

\pf. Let $X=B_{a_4r}\backslash B_{a_3r}$ and $d(x)$ be the distant
from $x\in X$ to $\mathbb{R}^n\backslash X$. By the elliptic
regularity, we obtain from \eqref{1.1} that $u\in
H^{2}_{loc}(\Omega\backslash \{0\})$. It is trivial that
\begin{equation}\label{3.2}
\begin{array}{l}
\| v\|_{H^{1}(\mathbb{R}^n)}\lesssim \|\Delta
v\|_{L^{2}(\mathbb{R}^n)} +\|v\|_{L^{2}(\mathbb{R}^n)}
\end{array}
\end{equation}
for all $v\in H^{2}(\mathbb{R}^n)$. By changing variables $x\to
B^{-1}x$ in \eqref{3.2}, we will have
\begin{equation}\label{3.3}
\begin{array}{l}
\sum_{|\alpha|\leq 1}B^{2-|\alpha|}\|
D^{\alpha}v\|_{L^{2}(\mathbb{R}^n)}\lesssim(\|\Delta
v\|_{L^{2}(\mathbb{R}^n)} +B^2\|v\|_{L^{2}(\mathbb{R}^n)})
\end{array}
\end{equation}
for all $v\in H^{2}(\mathbb{R}^n)$. To apply \eqref{3.3} on $u$, we
need to cut-off $u$. So let $\xi(x)\in C^{\infty}_0 ({\mathbb R}^n)$
satisfy $0\le\xi(x)\leq 1$ and
\begin{equation*}
\xi (x)=
\begin{cases}
\begin{array}{l}
1,\quad |x|<1/4,\\
0,\quad |x|\geq 1/2.
\end{array}
\end{cases}
\end{equation*}
Let us denote $\xi_y(x)=\xi((x-y)/d(y))$.   For $y\in X$, we apply
\eqref{3.3} to $\xi_y(x) u(x)$  and use \eqref{1.1} to get that
\begin{eqnarray}\label{3.4}
&& B^{2}\int_{|x-y|\leq d(y)/4}|\nabla u|^2dx\notag\\
&\lesssim&\int_{|x-y|\leq d(y)/2}(|A|^{2}+d(y)^{-2})|\nabla u|^2dx+
\int_{|x-y|\leq d(y)/2}|\nabla p|^2dx\notag\\
&& +(B^4+d(y)^{-4})\int_{|x-y|\leq d(y)/2}|u|^2dx.
\end{eqnarray}
Now taking $B=Md(y)^{-1}$ for some positive constant $M$ and
multiplying $d(y)^{4}$ on both sides of \eqref{3.4}, we have
\begin{eqnarray}\label{3.5}
&&M^{2}d(y)^{2}\int_{|x-y|\leq d(y)/4}|\nabla u|^2dx\notag\\
&\lesssim& \int_{|x-y|\leq
d(y)/2}(d(y)^{4}|A|^{2}+d(y)^{2})|\nabla u|^2dx\notag\\
&&+ \int_{|x-y|\leq d(y)/2}d(y)^{4}|\nabla
p|^2dx+(M^4+1)\int_{|x-y|\leq d(y)/2}|u|^2dx.
\end{eqnarray}

Integrating $d(y)^{-n}dy$ over $X$ on both sides of \eqref{3.5} and
using Fubini's Theorem, we get that
\begin{eqnarray}\label{3.6}
&&M^{2}\int_{X}\int_{|x-y|\leq d(y)/4}d(y)^{2-n}|\nabla u|^2dydx\notag\\
&\lesssim&\int_{X}\int_{|x-y|\leq
d(y)/2}(d(y)^{2}+d(y)^{4}|A|^{2})|\nabla
u(x)|^2d(y)^{-n}dydx\notag\\
&&\quad+ \int_{X}\int_{|x-y|\leq d(y)/2}d(y)^{4-n}|\nabla
p|^2dydx\notag\\
&&\quad+M^4\int_{X}\int_{|x-y|\leq d(y)/2}|u|^2d(y)^{-n}dydx.
\end{eqnarray}
Note that $|d(x)-d(y)|\leq |x-y|$. If $ |x-y|\leq d(x)/3$, then
\begin{equation}\label{3.7}
\begin{array}{l}
2d(x)/3\leq d(y)\leq 4d(x)/3.
\end{array}
\end{equation}
On the other hand, if $ |x-y|\leq d(y)/2$, then
\begin{equation}\label{3.8}
\begin{array}{l}
d(x)/2\leq d(y)\leq 3d(x)/2.
\end{array}
\end{equation}
By \eqref{3.7} and  \eqref{3.8}, we have
\begin{equation}\label{3.9}
\begin{cases}
\begin{array}{l}
\int_{|x-y|\leq d(y)/4}d(y)^{-n}dy\geq (3/4)^n\int_{|x-y|\leq
d(x)/6}d(x)^{-n}dy\geq 8^{-n}\int_{|y|\leq 1}dy,\\
\int_{|x-y|\leq d(y)/2}d(y)^{-n}dy\leq 2^n\int_{|x-y|\leq
3d(x)/4}d(x)^{-n}dy\leq (3/2)^{n}\int_{|y|\leq 1}dy
\end{array}
\end{cases}
\end{equation}
Combining \eqref{3.6}--\eqref{3.9}, we obtain
\begin{eqnarray}\label{3.10}
&&M^{2}\int_{X}d(x)^{2}|\nabla u|^2dx\notag\\
&\lesssim&\int_{X}(d(x)^{2}+d(x)^{4}|A|^{2})|\nabla u(x)|^2dx+
\int_{X}d(x)^{4}|\nabla
p|^2dx\notag\\
&&\hspace{5mm}+M^4\int_{X}|u|^2dx.
\end{eqnarray}
In view of \eqref{asing}, we can take $M$ large enough to absorb the
first term on the right hand side of \eqref{3.10}. Thus we conclude
that
\begin{equation}\label{3.11}
\int_{X}d(x)^{2}|\nabla u|^2dx \lesssim\int_{X}(d(x)^{4}|\nabla
p|^2+|u|^2)dx.
\end{equation}
We recall that $X=B_{a_4r}\backslash B_{a_3r}$ and note that
$d(x)\geq \tilde Cr$ if $x\in B_{a_2r}\backslash B_{a_1r}$, where
$\tilde C$ is independent of $r$. Hence, \eqref{3.1} is an easy
consequence of \eqref{3.11}.\eproof

\section{Proof of Theorem \ref{thm1.1} and Theorem \ref{thm1.2}}\label{sec4}
\setcounter{equation}{0}

This section is devoted to the proofs of Theorem \ref{thm1.1} and
Theorem \ref{thm1.2}. To begin, we first consider the case where
$0<R_1<R_2<R<1$ and $B_R\subset\Omega$. The small constant
$R$ will be determined later. Since $(u,p)\in (H^1(B_{R_0}))^{n+1}$,
the elliptic regularity theorem implies $u\in
H^2_{loc}(B_{R_0}\setminus \{0\})$. Therefore, to use estimate
\eqref{2.4}, we simply cut-off $u$. So let $\chi(x)\in C^{\infty}_0
({\mathbb R}^n)$ satisfy $0\le\chi(x)\leq 1$ and
\begin{equation*}
\chi (x)=
\begin{cases}
\begin{array}{l}
0,\quad |x|\leq R_1/e,\\
1,\quad R_1/2<|x|<eR_2,\\
0,\quad |x|\geq 3R_2,
\end{array}
\end{cases}
\end{equation*}
where $e=\exp(1)$. We remark that we first choose a small $R$ such
that $R\le\min\{r_0,r_1\}/3=\tilde R_0$, where $r_0$ and $r_1$ are
constants appeared in \eqref{2.4} and \eqref{2.14}. Hence $\tilde
R_0$ depends on $n$. It is easy to see that for any multiindex
$\alpha$
\begin{equation}\label{4.1}
\begin{cases}
|D^{\alpha}\chi|=O(R_1^{-|\alpha|})\ \text{for all}\ R_1/e\le |x|\le R_1/2\\
|D^{\alpha}\chi|=O(R_2^{-|\alpha|})\ \text{for all}\ eR_2\le |x|\le
3R_2.
\end{cases}
\end{equation}
Applying \eqref{2.4} to $\chi u$ gives
\begin{equation}\label{4.2}
{C}_1\beta  \int (\log|x|)^{-2}\varphi^2_\beta
|x|^{-n}(|x|^{2}|\nabla (\chi u)|^2+|\chi u|^2)dx \leq \int
\varphi^2_\beta |x|^{-n}|x|^{4}|\Delta (\chi u)|^2 dx.
\end{equation}
From now on, $C_1,C_2,\cdots$ denote general constants whose
dependence will be specified whenever necessary. Next applying
\eqref{2.14} to $w=\chi p$ and $f=|x|\chi A\cdot \nabla u$, we get
that
\begin{eqnarray}\label{2.19}
&&{C}_2\int\varphi^2_\beta (\log|x|)^2(|x|^{4-n}|\nabla (\chi p)|^2+|x|^{2-n}|\chi p|^2)dx\notag\\
&\leq& \beta\int \varphi^2_\beta
(\log|x|)^{4}|x|^{2-n}[|x|^{2}\Delta
(\chi p)+|x|{\rm div} (|x|\chi A\cdot \nabla u)]^2dx\notag\\
&&+\beta\int \varphi^2_\beta (\log|x|)^{4}|x|^{2-n}\| |x|\chi A\cdot
\nabla u\|^2 dx.
\end{eqnarray}
Multiplying by $M_1$ on \eqref{4.2} and combining \eqref{2.19}, we
obtain that
\begin{eqnarray}\label{4.3}
&&M_1\beta \int_{R_1/2<|x|<eR_2}
(\log|x|)^{-2}\varphi^2_\beta |x|^{-n}(|x|^{2}|\nabla u|^2+| u|^2)dx\notag\\
&&+\int_{R_1/2<|x|<eR_2}(\log|x|)^2\varphi^2_\beta |x|^{-n}(|x|^{4}|\nabla p|^2+|x|^{2}|p|^2)dx\notag\\
&\leq&M_1\beta \int
\varphi^2_\beta (\log|x|)^{-2}|x|^{-n}(|x|^{2}\nabla (\chi u)|^2+|\chi u|^2)dx\notag\\
&&+ \int(\log|x|)^2\varphi^2_\beta |x|^{-n}(|x|^{4}|\nabla (\chi p)|^2+|x|^{2}|\chi p|^2)dx\notag\\
&\leq &M_1C_{3}\int\varphi^2_\beta |x|^{-n}|x|^{4}|\Delta (\chi u)|^2 dx\notag\\
&&+\beta C_{3}\int(\log|x|)^{4}\varphi^2_\beta |x|^{-n}[|x|^{3}\Delta (\chi p)+|x|^2{\rm div} (|x|\chi A\cdot \nabla u)]^2dx\notag\\
&&+\beta C_{3}\int(\log|x|)^{4}\varphi^2_\beta |x|^{-n}\| |x|^2\chi
A\cdot \nabla u\|^2 dx.
\end{eqnarray}
By \eqref{1.1}, \eqref{asing}, \eqref{eq0}, and estimates
\eqref{4.1}, we deduce from \eqref{4.3} that
\begin{eqnarray}\label{4.4}
&&M_1\beta \int_{R_1/2<|x|<eR_2}
(\log|x|)^{-2}\varphi^2_\beta |x|^{-n}(|x|^{2}|\nabla u|^2+|u|^2)dx\notag\\
&&+\int_{R_1/2<|x|<eR_2}(\log|x|)^{2}\varphi^2_\beta |x|^{-n}(|x|^{4}|\nabla p|^2+|x|^{2}|p|^2)dx\notag\\
&\leq& C_4M_1\int_{R_1/2<|x|<eR_2}(\log|x|)^{-2}\varphi^2_\beta |x|^{-n}|x|^{2}|\nabla u|^2 dx\notag\\
&&+C_4M_1\int_{R_1/2<|x|<eR_2}\varphi^2_\beta |x|^{-n}|x|^{4}|\nabla p|^2 dx\notag\\
&&+C_{4}\beta \int_{R_1/2<|x|<eR_2}(\log|x|)^{-2}\varphi^2_\beta |x|^{-n}|x|^2|\nabla u|^2dx\notag\\
&&+C_4M_1\int_{\{R_1/e\le |x|\le R_1/2\}\cup\{eR_2\le |x|\le 3R_2\}}\varphi^2_\beta |x|^{-n}|\tilde{U}|^2 dx\notag\\
&&+C_{4}\beta \int_{\{R_1/e\le |x|\le R_1/2\}\cup\{eR_2\le |x|\le
3R_2\}}(\log|x|)^{4}\varphi^2_\beta |x|^{-n}|\tilde{U}|^2 dx,
\end{eqnarray}
where $|\tilde{U}(x)|^2=|x|^{4}|\nabla
p|^2+|x|^{2}|p|^2+|x|^{2}|\nabla u|^2+|u|^2$ and the positive
constant $C_4$ only depends on $n$.

Now letting $M_1=2+2C_4$, $\beta\geq 2+2C_{4}$, and $R$ small enough
such that $(\log(eR))^2\ge 2C_4M_1$, then the first three terms on
the right hand side of \eqref{4.4} can be absorbed by the left hand
side of \eqref{4.4}. Also, it is easy to check that there exists
$\tilde R_1>0$, depending on $n$, such that for all $\beta>0$, both
$(\log|x|)^{-2}|x|^{-n}\varphi_{\beta}^2(|x|)$ and
$(\log|x|)^{4}|x|^{-n}\varphi_{\beta}^2(|x|)$ are decreasing
functions in $0<|x|<\tilde R_1$. So we choose a small $R<\tilde
R_2$, where $\tilde R_2=\min\{\exp(-2\sqrt{2C_4M_1}-1),\tilde
R_1/3,\tilde R_0\}$. It is clear that $\tilde R_2$ depends on $n$.
With the choices described above, we obtain from \eqref{4.4} that

\begin{eqnarray}\label{4.5}
&&R_2^{-n}(\log
R_2)^{-2}\varphi^2_\beta(R_2)\int_{R_1/2<|x|<R_2}|\tilde{U}|^2dx\notag\\
&\leq &\int_{R_1/2<|x|<eR_2}(\log|x|)^{-2}\varphi^2_\beta |x|^{-n}|\tilde{U}|^2dx\notag\\
&\leq &C_5\beta \int_{\{R_1/e\le |x|\le R_1/2\}\cup\{eR_2\le |x|\le 3R_2\}}(\log|x|)^{4}\varphi^2_\beta |x|^{-n}|\tilde{U}|^2 dx\notag\\
&\leq &C_5\beta(\log(R_1/e))^{4}(R_1/e)^{-n}\varphi^2_\beta(R_1/e)\int_{\{R_1/e\le |x|\le R_1/2\}}|\tilde{U}|^2 dx\notag\\
&&+C_5\beta(\log(eR_2))^{4}(eR_2)^{-n}\varphi^2_\beta(eR_2)\int_{\{eR_2\le
|x|\le 3R_2\}}|\tilde{U}|^2 dx.
\end{eqnarray}

Using \eqref{3.1}, we can control $|\nabla u|$ terms on the right
hand side of \eqref{4.5}. In other words, it follows from
\eqref{3.1} that
\begin{eqnarray}\label{4.6}
&&R_2^{-2\beta-n}(\log R_2)^{-4\beta-2}\int_{R_1/2<|x|<R_2}|U|^2dx\notag\\
&\leq & C_{6}2^{2\beta+n}(\log(R_1/e))^{4}(R_1/e)^{-n}\varphi^2_\beta(R_1/e)\int_{\{R_1/4\le |x|\le R_1\}}|U|^2 dx\notag\\
&&+C_{6}2^{2\beta+n} (\log(eR_2))^{4}(eR_2)^{-n}\varphi^2_\beta(eR_2)\int_{\{2R_2\le |x|\le 4R_2\}}|U|^2 dx\notag\\
&= & C_{6}2^{2\beta+n}(\log(R_1/e))^{-4\beta+4}(R_1/e)^{-2\beta-n}\int_{\{R_1/4\le |x|\le R_1\}}|U|^2 dx\notag\\
&&+C_{6}2^{2\beta+n}
(\log(eR_2))^{-4\beta+4}(eR_2)^{-2\beta-n}\int_{\{2R_2\le |x|\le
4R_2\}}|U|^2 dx.
\end{eqnarray}
Recall that $|U(x)|^2=|x|^{4}|\nabla p|^2+|x|^{2}|p|^2+|u|^2$.
Replacing $2\beta+n$ by $\beta$, \eqref{4.6} becomes
\begin{eqnarray}\label{4.7}
&&R_2^{-\beta}(\log R_2)^{-2\beta+2n-2}\int_{R_1/2<|x|<R_2}|U|^2dx\notag\\
&\leq & C_{7}2^\beta(\log(R_1/e))^{-2\beta+2n+4}(R_1/e)^{-\beta}\int_{\{R_1/4\le |x|\le R_1\}}|U|^2 dx\notag\\
&&+C_{7}2^\beta
(\log(eR_2))^{-2\beta+2n+4}(eR_2)^{-\beta}\int_{\{2R_2\le |x|\le
4R_2\}}|U|^2 dx.
\end{eqnarray}
Dividing $R_2^{-\beta}(\log R_2)^{-2\beta+2n-2}$ on the both sides
of \eqref{4.7} and providing $\beta\geq n+2$, we have that
\begin{eqnarray}\label{4.8}
&&\int_{R_1/2<|x|<R_2}|U|^2dx\notag\\
&\leq & C_{8}(\log R_2)^6(2eR_2/R_1)^{\beta}\int_{\{R_1/4\le |x|\le R_1\}}|U|^2 dx\notag\\
&&+C_{8}(\log R_2)^{6}(2/e)^\beta [(\log R_2/\log(eR_2))^{2}]^{\beta-n-2}\int_{\{2R_2\le |x|\le 4R_2\}}|U|^2 dx\notag\\
&\leq & C_{8}(\log R_2)^{6}(2eR_2/R_1)^{\beta}\int_{\{R_1/4\le |x|\le R_1\}}|U|^2 dx\notag\\
&&+C_{8}(\log R_2)^{6}(4/5)^\beta \int_{\{2R_2\le |x|\le
4R_2\}}|U|^2 dx.
\end{eqnarray}
In deriving the second inequality above, we use the fact that
$$
\frac{\log R_2}{\log(eR_2)}\to 1\quad\text{as}\quad R_2\to 0,
$$
and thus
$$
\frac{2}{e}\cdot\frac{\log R_2}{\log(eR_2)}<\frac{4}{5}
$$
for all $R_2<\tilde R_3$, where $\tilde R_3$ is sufficiently small.
We now take $\tilde R=\min\{\tilde R_2,\tilde R_3\}$, which depends
on $n$.

Adding $\int_{|x|<{R_1/2}} |U|^2 dx$ to both sides of \eqref{4.8}
leads to
\begin{eqnarray}\label{4.9}
\int_{|x|<R_2}|U|^2dx&\leq&  C_{9}(\log
R_2)^{6}(2eR_2/R_1)^{\beta}\int_{|x|\le R_1}|U|^2
dx\notag\\
&&+C_{9}(\log R_2)^{6}(4/5)^\beta \int_{|x|\le 1}|U|^2 dx.
\end{eqnarray}
It should be noted that \eqref{4.9} holds for all
$\beta\ge\tilde\beta$ with $\tilde\beta$ depending only on $n$. For
simplicity, by denoting
\begin{equation*}
E(R_1,R_2)=\log(2eR_2/R_1),\quad B=\log (5/4),
\end{equation*}
\eqref{4.9} becomes
\begin{eqnarray}\label{4.10}
&&\int_{|x|<R_2}|U|^2dx\notag\\
&\leq & C_{9}(\log R_2)^{6}\Big{\{}\exp(E\beta)\int_{|x|<{R_1}}
|U|^2 dx+\exp(-B\beta) \int_{|x|<1} |U|^2 dx\Big{\}}.\notag\\
\end{eqnarray}

To further simplify the terms on the right hand side of
\eqref{4.10}, we consider two cases. If $\int_{|x|<{R_1}} |U|^2 dx\ne 0$ and
$$\exp{(E\tilde\beta)}\int_{|x|<{R_1}} |U|^2 dx<\exp{(-B\tilde\beta)}\int_{|x|<{1}} |U|^2 dx,$$
then we can pick a $\beta>\tilde\beta$ such that
$$
\exp{(E\beta)}\int_{|x|<{R_1}} |U|^2
dx=\exp{(-B\beta)}\int_{|x|<{1}} |U|^2 dx.
$$
Using such $\beta$, we obtain from \eqref{4.10} that
\begin{eqnarray}\label{4.11}
&&\int_{|x|<R_2}|U|^2dx\notag\\
&\leq& 2C_{9}(\log R_2)^{6}\exp{(E\beta)} \int_{|x|<{R_1}} |U|^2dx\notag\\
&=& 2C_{9}(\log
R_2)^{6}\left(\int_{|x|<{R_1}}|U|^2dx\right)^{\frac{B}{E+B}}\left(\int_{|x|<{1}}|U|^2dx\right)^{\frac{E}{E+B}}.
\end{eqnarray}
If $\int_{|x|<{R_1}} |U|^2 dx= 0$, then letting $\beta\to\infty$ in \eqref{4.10} we have $\int_{|x|<R_2}|U|^2dx=0$ as well. The three-ball inequality obviously holds.

On the other hand, if
$$ \exp{(-B\tilde\beta)}\int_{|x|<{1}} |U|^2dx\leq\exp{(E\tilde\beta)}\int_{|x|<{R_1}} |U|^2 dx,$$
then we have
\begin{eqnarray}\label{4.12}
&&\int_{|x|<{R_2}}|U|^2 dx\notag\\
&\leq & \left(\int_{|x|<1}|U|^2dx\right)^{\frac{B}{E+B}}\left(\int_{|x|<1}|U|^2dx\right)^{\frac{E}{E+B}}\notag\\
&\leq &
\exp{(B\tilde\beta)}\left(\int_{|x|<{R_1}}|U|^2dx\right)^{\frac{B}{E+B}}\left(\int_{|x|<1}|U|^2dx\right)^{\frac{E}{E+B}}.
\end{eqnarray}
Putting together \eqref{4.11}, \eqref{4.12}, and setting
$C_{10}=\max\{2C_{9}(\log R_2)^{6},\exp{(\tilde\beta\log(5/4))}\}$, we arrive at
\begin{equation}\label{4.13}
\int_{|x|<{R_2}}|U|^2 dx \le
C_{10}\left(\int_{|x|<{R_1}}|U|^2dx\right)^{\frac{B}{E+B}}\left(\int_{|x|<1}|U|^2dx\right)^{\frac{E}{E+B}}.
\end{equation}
It is readily seen that $\frac{B}{E+B}\approx (\log
(1/R_1))^{-1}$ when $R_1$ tends to $0$.

Now for the general case, we consider $0<R_1<R_2<R_3<1$ with
$R_1/R_3<R_2/R_3\le \tilde R$, where $\tilde R$ is given as above. By scaling, i.e. defining
$\widehat{u}(y):=u(R_3y)$, $\widehat{p}(y):=R_3p(R_3y)$ and
$\widehat{A}(y)=A(R_3y)$,  \eqref{4.13}  becomes
\begin{equation}\label{4.14}
\int_{|y|<{R_2/R_3}}|\widehat{U}(y)|^2 dy \leq
C_{11}(\int_{|y|<{R_1/R_3}}|\widehat{U}(y)|^2dy)^{\tau}(\int_{|y|<1}|\widehat{U}(y)|^2dy)^{1-\tau},
\end{equation}
where $$\tau=B/[E(R_1/R_3,R_2/R_3)+B],$$ $$C_{11}=\max\{2C_{9}(\log
R_2/R_3)^{6},\exp{(\tilde\beta\log(5/4))}\},$$ and $\widehat{U}(y)=|y|^{4}|\nabla
\widehat{p}(y)|^2+|y|^{2}|\widehat{p}(y)|^2+|\widehat{u}(y)|^2$. Note that $C_{11}$ is independent of $R_1$. Restoring the
variable $x=R_3y$ in \eqref{4.14} gives
\begin{equation*}
\int_{|x|<{R_2}}|U|^2 dx \leq
C_{11}(\int_{|x|<{R_1}}|U|^2dx)^{\tau}(\int_{|x|<{R_3}}|U|^2dx)^{1-\tau}.
\end{equation*}
The proof of Theorem \ref{thm1.1} is complete.

We now turn to the proof of Theorem \ref{thm1.2}. We fix $R_2$, $R_3$ in Theorem~\ref{thm1.1} and  define
$$
\begin{cases}
\widetilde{u}(x):=u(x)/\sqrt{\int_{|x|<{R_2}}|U|^2 dx},\\
\widetilde{p}(x):=p(x)/\sqrt{\int_{|x|<{R_2}}|U|^2 dx},\\
 V(x)=|x|^{4}|\nabla
\widetilde{p}(x)|^2+|x|^{2}|\widetilde{p}(x)|^2+|\widetilde{u}(x)|^2.
\end{cases}
$$
Note that $\int_{|x|<{R_2}}|V|^2 dx=1$. From the three-ball inequality
\eqref{1.2}, we have that
\begin{equation}\label{4.16}
1 \leq
C(\int_{|x|<{R_1}}|V|^2dx)^{\tau}(\int_{|x|<{R_3}}|V|^2dx)^{1-\tau}.
\end{equation}
Raising both sides by $1/\tau$ yields that
\begin{equation}\label{4.17}
\int_{|x|<{R_3}}|V|^2dx \leq (\int_{|x|<{R_1}}|V|^2dx)(C\int_{|x|<{R_3}}|V|^2dx)^{1/\tau}.
\end{equation}
In view of the formula for $\tau$, we can deduce from \eqref{4.17} that
\begin{equation}\label{4.18}
\int_{|x|<{R_3}}|V|^2dx \leq (\int_{|x|<{R_1}}|V|^2dx)(1/R_1)^{\tilde C\log(\int_{|x|<{R_3}}|V|^2dx)},
\end{equation}
where $\tilde C$ is a positive constant depending on $n$ and
$R_2/R_3$. Consequently, \eqref{4.18} is equivalent to
$$
(\int_{|x|<R_3}|U|^2dx) R_1^m\le \int_{|x|<R_1}|U|^2dx
$$
for all $R_1$ sufficiently small, where
$$
m=\tilde C\log\Big{(}\frac{\int_{|x|<R_3}|U|^2dx}{\int_{|x|<R_2}|U|^2dx}\Big{)}.
$$
We now end the proof of Theorem \ref{thm1.2}.

Finally, we come to the proof of Corollary~\ref{cor1}. In view of
Theorem~\ref{thm1.2}, it is enough to show that
\begin{equation}\label{4.19}
p\in H^1_{loc}(\Omega)
\end{equation}
and for all $N>0$
\begin{equation}\label{4.20}
\int_{|x|<R}|\nabla p|^2dx=O(R^N)\quad\text{as}\quad R\to 0.
\end{equation}
It is only a matter of checking that the arguments used in
\cite[page 1898-1899]{Reg2} can be applied to prove \eqref{4.19} and
\eqref{4.20}. To avoid unnecessary repetition, we only sketch the
main steps here. By virtue of \eqref{eq0}, $u\in H^1_{loc}(\Omega)$,
and \eqref{asing}, we get that $p\in
H^1_{loc}(\Omega\setminus\{0\})$ by elliptic regularity. Using
elliptic regularity and the first equation of \eqref{1.1}, we have
$u\in H^2_{loc}(\Omega\setminus\{0\})$. By the vanishing assumption
\eqref{info}, we can derive that
\begin{equation*}
\int_{R<|x|<2R}|\nabla p|^2dx=O(R^N)\quad\text{as}\quad R\to 0
\end{equation*}
for all $N>0$. It follows that $p$ is the sum of a function in
$H^1_{loc}(\Omega)$ and a distribution supported at $0$. But no
distribution supported at $0$ is in $L^2_{loc}(\Omega)$. Thus, $p\in
H^1_{loc}(\Omega)$ and \eqref{4.20} holds.

\section*{Acknowledgements}
The authors were supported in part by the National Science Council
of Taiwan.


\begin{thebibliography}{50}
\bibitem{abrv}
 G. Alessandrini, E. Beretta, E. Rosset, and S. Vessella, \emph{Optimal stability for elliptic
 boundary value problems with unknow boundaries}, Ann. Scuola Norm. Sup. Pisa Cl. Sci, {\bf 29} (2000), 755-786.

\bibitem{bou}
J. Bourgain and C. Kenig, \emph{On localization in the continuous
Anderson-Bernoulli model in higher dimension}, Invent. Math. {\bf
161} (2005), 1389-426.

\bibitem{dofe}
H. Donnelly and C. Fefferman, \emph{Nodal sets of eigenfunctions on
Riemannian manifolds}, Invent. Math. {\bf 93} (1988), 161-183.

\bibitem{efv}
L. Escauriaza, F.J. Fern\'andez, and S. Vessella, \emph{Doubling
properties of caloric functions}, Appl. Anal., \textbf{85} (2006),
205-223.
\bibitem{fale}
C. Fabre and G. Lebeau, \emph{Prolongement unique des solutions de
l'\'{e}quation de Stokes}, Comm. in PDE, \textbf{21} (1996),
573-596.
\bibitem{gl1}
N. Garofalo and F.H. Lin, \emph{Monotonicity properties of
variational integrals, $A_p$ weights and unique continuation},
Indiana Univ. Math. J. {\bf 35} (1986), 245-267.
\bibitem{gl2}
N. Garofalo and F.H. Lin,  \emph{Unique continuation for elliptic
operators: a geometric-variational approach}, Comm. Pure Appl.
Math., {\bf 40}, 347-366, 1987.
\bibitem{hasi}
R. Hardt and L. Simon, \emph{Nodal sets for solutions of elliptic
equations}, J. Diff. Geom., \textbf{30} (1989), 505-522.
\bibitem{Hor1}
L. H\"{o}rmander, {\it Uniqueness theorems for second order elliptic
differential equations}, Comm. in P.D.E. {\bf 8}, No. 1, (1983),
21-64.
\bibitem{Hor3}
L. H\"{o}rmander, {"The analysis of linear partial differential
operators"}, Vol. 3, Springer-Verlag, Berlin/New York, 1985.
\bibitem{jeke}
D. Jerison and C. Kenig, \emph{Unique continuation and absence of
positive eigenvalues for Schrodinger operators. With an appendix by
E. M. Stein}, Ann. of Math. (2), \textbf{121} (1985), 463-494.
\bibitem{kenigcna}
C. Kenig, \emph{Lecture Notes for 2006 CNA Summer School:
Probabilistic and Analytical Perspectives on Contemporary PDEs},
Center for Nonlinear Analysis, Carnegie Mellon University.
\bibitem{lin1}
F.H. Lin, \emph{A uniqueness theorem for parabolic equations}, Comm.
Pure Appl. Math., \textbf{43} (1990), 127-136.
\bibitem{linf2}
F.H. Lin, \emph{Nodal sets of solutions of elliptic and parabolic
equations}, Comm. Pure Appl. Math., \textbf{44} (1991), 287-308.
\bibitem{lin2}
C.L. Lin, G. Nakamura and J.N. Wang \emph{Quantitative uniqueness
for second order elliptic operators with strongly singular
coefficients}, Preprint.
\bibitem{lin3}
C.L. Lin, S. Nagayasu and J.N. Wang  \emph{Quantitative uniqueness
for the power of Laplacian with singular coefficients}, Preprint.
\bibitem{Reg}
R. Regbaoui, {\it Strong uniqueness for second order differential
operators}, J. Diff. Eq. {\bf 141} (1997), 201--217.
\bibitem{Reg2}
R. Regbaoui, {\it Strong unique continuation for stokes equations},
Comm. in PDE {\bf 24} (1999), 1891-1902.

\end{thebibliography}
\end{document}